\documentclass[11pt]{article}

\def\ArticleTitle{Benchmarking Bilevel Derivative-Free Optimization Algorithms}
\def\FirstAuthorName{Valentin Dijon}
\def\FirstAuthorEmail{valentin.dijon@etud.polymtl.ca}
\def\FirstAuthorInstitution{GERAD and Department of Mathematics and Industrial Engineering, Polytechnique Montr\'eal}

\def\SecondAuthorName{Youssef Diouane}		
\def\SecondAuthorEmail{youssef.diouane@polymtl.ca}
\def\SecondAuthorInstitution{GERAD and Department of Mathematics and Industrial Engineering, Polytechnique Montr\'eal}

\def\ThirdAuthorName{Charles Audet}		
\def\ThirdAuthorEmail{charles.audet@polymtl.ca}
\def\ThirdAuthorInstitution{GERAD and Department of Mathematics and Industrial Engineering, Polytechnique Montr\'eal}

\usepackage{graphicx}%
\graphicspath{{Figures/}} 
\usepackage{multirow}%
\usepackage{amsmath,amssymb,amsfonts}%
\usepackage{amsthm}
\usepackage{mathrsfs}%
\usepackage[title]{appendix}%
\usepackage{xcolor}%
\usepackage{textcomp}%
\usepackage{manyfoot}%
\usepackage{booktabs}%
\usepackage{algorithm}%
\usepackage{algorithmicx}%
\usepackage{algpseudocode}%
\usepackage{listings}%

\usepackage[margin=1in]{geometry}
\usepackage{enumerate}
\usepackage{mathtools}
\usepackage{subdepth}  
\usepackage{subcaption}
\usepackage{array}
\usepackage{dsfont}
\usepackage{hyperref} 
\usepackage[capitalise]{cleveref}
\usepackage{longtable}
\usepackage[textwidth=22mm, backgroundcolor=yellow, linecolor=orange, textsize=scriptsize]{todonotes}
\usepackage[short,c2]{optidef}
\DeclareMathOperator{\st}{s.t.}

\usepackage{pifont}
\usepackage{longtable} 
\usepackage[normalem]{ulem}

\usepackage{tikz}
\usepackage{standalone}

\usepackage{pgfplots}
\pgfplotsset{compat=newest}
\pgfplotsset{every axis legend/.append style={%
cells={anchor=west}}
}
\usetikzlibrary{arrows}
\tikzset{>=stealth'}

\pgfplotsset{
  seriesBlank/.style={color=white,  thick, solid, line width=0.5pt},
  seriesA/.style={color=blue,  thick, solid, line width=0.5pt,
                  mark=triangle},
  seriesB/.style={color=red,   thick, solid, line width=0.5pt,
                  mark=square},
  seriesC/.style={color=green, thick, solid, line width=0.5pt,
                  mark=diamond},
  seriesD/.style={color=orange, thick, solid, line width=0.5pt,
                  mark=o},
  seriesE/.style={color=purple, thick, solid, line width=0.5pt,
                  mark=x},
}

\newcommand{\globallegendInternReferee}{%
  \begin{tikzpicture}
    \begin{axis}[
      hide axis,
      xmin=0, xmax=1, ymin=0, ymax=1,
      width=16cm, height=10cm,       
      legend columns=5,
      legend style={
        draw=black,
        fill=white,
        inner sep=6pt,
        column sep=18pt,
        font=\small,
        /tikz/every even column/.append style={column sep=6pt},
      },
      legend to name={globallegendInternReferee}, 
    ]
      \addlegendimage{seriesBlank}\addlegendentry{\hspace{-1.2cm} LL subsolver:}
      \addlegendimage{seriesA}\addlegendentry{\hspace{-0.5cm} \texttt{Algo1}}
      \addlegendimage{seriesB}\addlegendentry{\hspace{-0.5cm} \texttt{Algo2}}
      \addlegendimage{seriesC}\addlegendentry{\hspace{-0.5cm} \texttt{Algo3}}
    \end{axis}
  \end{tikzpicture}%
}

\newcommand{\globallegendExternReferee}{%
  \begin{tikzpicture}
    \begin{axis}[
      hide axis,
      xmin=0, xmax=1, ymin=0, ymax=1,
      width=16cm, height=10cm,       
      legend columns=5,
      legend style={
        draw=black,
        fill=white,
        inner sep=6pt,
        column sep=18pt,
        font=\small,
        /tikz/every even column/.append style={column sep=6pt},
      },
      legend to name={globallegendExternReferee}, 
    ]
      \addlegendimage{seriesA}\addlegendentry{\hspace{-0.5cm} \texttt{Algo1}}
      \addlegendimage{seriesB}\addlegendentry{\hspace{-0.5cm} \texttt{Algo2}}
      \addlegendimage{seriesC}\addlegendentry{\hspace{-0.5cm} \texttt{Algo3}}
    \end{axis}
  \end{tikzpicture}%
}

\newcommand{\globallegendCompareLambda}{%
  \begin{tikzpicture}
    \begin{axis}[
      hide axis,
      xmin=0, xmax=1, ymin=0, ymax=1,
      width=16cm, height=10cm,       
      legend columns=5,
      legend style={
        draw=black,
        fill=white,
        inner sep=6pt,
        column sep=18pt,
        font=\small,
        /tikz/every even column/.append style={column sep=6pt},
      },
      legend to name={globallegendCompareLambda}, 
    ]
      \addlegendimage{seriesA}\addlegendentry{\hspace{-0.5cm} \texttt{Algo1}}
      \addlegendimage{seriesC}\addlegendentry{\hspace{-0.5cm} \texttt{Algo2}}
    \end{axis}
  \end{tikzpicture}%
}

\newcommand*\OK{\ding{51}}

\renewcommand{\mid}{\,:\,}

\newtheorem{definition}{Definition}

\usepackage[
  createShortEnv,
  conf={one big link},
  commandRef=Cref]{proof-at-the-end}

\makeatletter
\def\thmt@innercounters{section,equation,theorem}
\makeatother

\newcommand{\intInterval}[2]{[\![ #1,#2 ]\!]}

\NewDocumentEnvironment{proofEE}{O{}+b}{%
  \begin{proofE}[#1]
    #2
    \space
  \end{proofE}
}{}
\usepackage[numbers,sort&compress]{natbib}
\makeatletter
\renewcommand\NAT@bibsetnum[1]{\settowidth\labelwidth{\@biblabel{#1}}%
   \setlength{\leftmargin}{\bibindent}\addtolength{\leftmargin}{\dimexpr\labelwidth+\labelsep\relax}%
   \setlength{\itemindent}{-\bibindent}%
   \setlength{\listparindent}{\itemindent}
\setlength{\itemsep}{\bibsep}\setlength{\parsep}{\z@}%
   \ifNAT@openbib
     \addtolength{\leftmargin}{\bibindent}%
     \setlength{\itemindent}{-\bibindent}%
     \setlength{\listparindent}{\itemindent}%
     \setlength{\parsep}{0pt}%
   \fi
}
\makeatother

\makeatletter
\renewcommand*{\NAT@spacechar}{~} 
\makeatother

\usepackage{etoolbox}
\makeatletter
\patchcmd{\NAT@test}{\else \NAT@nm}{\else \NAT@hyper@{\NAT@nm}}{}{}
\makeatother
\usepackage{dsfont}
\newcommand{\R}{\mathds{R}}
\newcommand{\N}{\mathds{N}}

\newcommand{\A}{\mathcal{A}}
\newcommand{\Pbank}{\mathcal{P}}
\newcommand{\IR}{\mathrm{IR}}
\newcommand{\Refs}{\mathcal{R}}

\newcommand{\argminim}{\mathop{\textup{argmin}}}
\newcommand{\argmin}[1]{\argminim_{#1}}

\title{\ArticleTitle}
\author{%
\ThirdAuthorName\footnote{%
\ThirdAuthorInstitution. E-mail: \href{mailto:\ThirdAuthorEmail}{\ThirdAuthorEmail}.
}
\and
\FirstAuthorName\footnote{%
\FirstAuthorInstitution. E-mail: \href{mailto:\FirstAuthorEmail}{\FirstAuthorEmail}.
  }
  \and 
  \SecondAuthorName\footnote{%
  \SecondAuthorInstitution. E-mail: \href{mailto:\SecondAuthorEmail}{\SecondAuthorEmail}.
  }
}

\newcommand{\rd}{\color{red}}
\newcommand{\bl}{\color{black}}

\begin{document}

\maketitle

\begin{abstract}

Bilevel optimization involves an upper-level and a lower-level decision maker.
The lower-level optimization problem is nested within the constraints of the upper-level one.
A point is said to be admissible for the bilevel problem if it satisfies all constraints and is optimal for the lower-level decision-maker.
Bilevel derivative-free optimization (BL-DFO) algorithms address bilevel optimization problems in which {\color{black} the upper-level is solved using a derivative-free optimization method, and an inexact method it used to solve the lower- level.}
In this context, existing BL-DFO benchmarking techniques often do not rigorously validate the admissibility of proposed solutions, and do not adequately account for the computational effort deployed by the upper- and lower-level solvers.
This work proposes a benchmarking methodology for BL-DFO algorithms. 
A post-optimization procedure, named the refereeing procedure, is introduced to discard non-admissible points and ensures a fair comparison between the algorithms.
The computational effort deployed by upper- and lower-level solvers is also taken into account into the overall computational cost.
Numerical experiments illustrate the benchmarking methodology.\\


\end{abstract}


       

\section{Introduction}

Bilevel optimization (BLO) is a hierarchical optimization framework in which a subset of the decision variables is controlled by a second decision-maker, guided by its own objective function and subject to its own constraints. 
BLO has been extensively studied in the literature (e.g.~\citep{HaJaSa91,Dempe2002,CoMaSa2007,KlLaLjSc2021}) and arises in various application domains. 
BLO naturally models the leader-follower interactions from game theory, especially Stackelberg games~\citep{VoSPeScHu1953,Papavassilopoulos1982,Maheshwari2023}.
Other applications in engineering also arise from this structure when the two problems are in conflict~\citep{AgMSAbMi2012,LeBo1986,ZhLi2014}.
In recent years, BLO has gained renewed attention due to its relevance in machine learning and artificial intelligence.
Indeed, problems such as hyper-parameter optimization~\citep{bao2021stability,Franceschi2018}, meta-learning~\citep{Franceschi2018}, generative adversarial networks~\citep{goodfellow2014generative,gulrajani2017improved}, or robust neural network training~\citep{Madry2017} are naturally formalized with a bilevel structure.

Denote by $x$ the upper-level variables and by $y$ the lower-level variables.
Denote by $n_x$ and $n_y$ the dimensions of $x$ and $y$ respectively.
A bilevel problem is formulated in its general (optimistic) form by
\begin{equation}
    \begin{aligned}
    \min_{x \in X, y^* \in \R^{n_y}} \quad  & F(x,y^*)  \\
    \st \quad & G(x,y^*)  \le 0 \\
      & y^* \in   \argmin {y \in Y} \quad  f(x,y)\\
      & \quad \quad \quad \; \st \quad \; g(x,y) \le 0, \\
    \end{aligned}
    \label{eq:MainBilevelProb_Optimistic}
\end{equation}
where $X\subseteq \R^{n_x}$ and $Y \subseteq \R^{n_y}$ represent the respective domains of $x$ and $y$.
The upper-level (UL) problem is represented by an objective function $F : X \times Y \to \R$ and constraints $G: X \times Y \to \R^{n_G}$, where $n_G$ is the dimension of the upper-level (coupling) constraints.
Conversely, the lower-level (LL) problem is described by its own objective function $f: X \times Y \to \R$ and constraints $g: X \times Y \to \R^{n_g}$, where $n_g$ is the dimension of the lower-level constraints.
{\bl The lower-level feasible set associated with an upper-level variable $x \in X$ is denoted by $Y(x) \coloneqq \{y \in Y \mid g(x,y) \le 0\}$.}

A BLO problem involves a hierarchical decision process where the upper-level decision-maker selects a variable $x \in X$. 
However, evaluating the upper-level objective $F$ and constraints $G$ requires anticipating the reaction of the lower-level decision-maker.
Given a fixed $x$, the lower-level problem requires to minimize $f(x,y)$ {\bl over $Y(x)$}.
The set of optimal solutions to this problem is not necessarily a singleton.
In such cases, the upper-level may select any $y$ among the lower-level optimal solutions for which the follower is indifferent.
Among all such lower-level optimal solutions, the upper-level decision-maker selects one that satisfies the constraints $G(x,y) \leq 0$ and minimizes the objective $F(x,y)$.
This setting is referred to as optimistic bilevel programming.

A standard approach to solve the BLO problem~\eqref{eq:MainBilevelProb_Optimistic} consists in reformulating it as a single-level problem either by replacing the lower-level problem with its KKT optimality conditions~\citep{Bard1983} or by embedding the hyper-function associated with the lower-level value function (LLVF)~\citep{Dempe1992,Outrata1990}.
Penalization-based first-order methods built upon such reformulations have been shown to be effective~\citep{ZeZh2021}.
Such methods, however, require first-order information or restrictive assumptions on the upper-level, which 
are not always available.

This work focuses on benchmarking techniques employed for comparing bilevel derivative-free optimization (BL-DFO) algorithms where the {\bl upper-level} problem in~\eqref{eq:MainBilevelProb_Optimistic} is solved using a derivative-free optimization (DFO) method.
The BL-DFO methods considered proceed iteratively by proposing an upper-level variable $x$ and computing {\bl inexactly} its associated optimal lower-level response $y$.
Each iteration then involves solving the lower-level problem for a given $x$.
A wide range of methods can be applied to solve the lower-level problem.
Ensuring the optimality of lower-level responses $y$ with respect to a given $x$ is challenging and must be accounted for by benchmarking techniques in BLO.

\subsection{Benchmarking techniques in BL-DFO}

Benchmarking methodologies in DFO are well established for single-level DFO, but their extension to BLO is not trivial~\citep[Chapter 4]{AuHa2026}.
A basic approach consists in reporting performance tables that provide raw data of optimization runs across solvers and test instances.
Such tabular comparisons are used to support the performance of BL-DFO methods, e.g.,~\cite{DoPr2023,Kieffer2017}.
However, relying exclusively on tables may hinder interpretability, as extracting consistent conclusions from extensive tabulated results is cumbersome.
The difficulty increases substantially when multiple solvers, runs, and test instances are considered.
Convergence plots provide a graphical representation of the information reported in tables of statistics, namely the evolution of the value of the objective function with respect to the number of evaluations.
They are widely used to compare BL-DFO algorithms (see e.g., \cite{CoVi2012} and references therein) to illustrate their performances on selected instances.
Such a representation is appropriate when analyzing the performance of a given algorithm on a specific problem.

Although performance~\citep{DoMo02} and data~\citep{MoWi2009} profiles are commonly used in DFO for benchmarking algorithms, their use remains limited in the BL-DFO literature.
\citet{Islam2017} assess the performance of their population-based DFO algorithm by multiple means, including performance profiles. 
In the context of direct search methods,~\citet{DiKuRiZe2024} use data profiles to compare both mesh-based and sufficient-decrease schemes, constructing the profiles solely using upper-level function evaluations.
Recently, the line-search method devised by~\citet{Cesaroni2026} is benchmarked with performance and data profiles using lower-level evaluations as budget units.

The aforementioned benchmarking techniques are designed for single-level DFO and do not account for issues specific to BLO.
{\bl First, existing BL-DFO methods enforce lower-level optimality through additional algorithmic mechanisms~\citep{Cesaroni2026,Islam2017}, but no unified benchmarking methodology exists to verify whether the solutions reported by a solver are in fact feasible and lower-level optimal.
Such verification is essential for a fair comparison between BL-DFO algorithms. Second, the measures of computational effort used to compare solvers do not always count both upper- and lower-level function evaluations.} 
The upper- and the lower-level efforts need to be adequately aggregated to represent the true computational cost deployed to solve~\eqref{eq:MainBilevelProb_Optimistic}.

\subsection{Motivation and contributions}

The main difficulty addressed in this paper arises when a DFO method claims that a point $(x,y)$ is admissible for a bilevel problem, i.e., that it is feasible and that $y$ is an optimal solution of the lower-level problem with respect to $x$.
{\bl Ensuring the admissibility of points is crucial to avoid misleading conclusions about the performance of BL-DFO algorithms. For instance, consider the bilevel problem
\begin{equation}
    \begin{aligned}
    \min_{x \in \R,\, y \in \R} \quad & F(x,y) = 2x - y \\
    \st \quad & x \ge 0, \\
              & y \in \argmin {z \in Y(x)} \ z,
    \end{aligned}
    \label{eq:ExampleBilevelProb_Optimistic}
\end{equation}
where $Y(x) = [x,5]$.
For any $x \in [0,5]$, the lower-level problem has the unique solution $y = x$, so the inducible region is $\mathrm{IR} = \{(x,x) : 0 \le x \le 5\}$, represented by the green line in~\Cref{fig:ExampleNeedAdmissibility}.
The set $\Omega = \{(x,y) : x \ge 0,\ y \in Y(x)\}$, with upper- and lower-level constraints, is represented in blue.
On IR, the upper-level objective reduces to $F(x,x) = x$, so the optimal solution is $(x^*,y^*) = (0,0)$, represented by the green square, with $F(0,0) = 0$.

\begin{figure}[ht]
    \centering
    \includegraphics[width=0.45\linewidth]{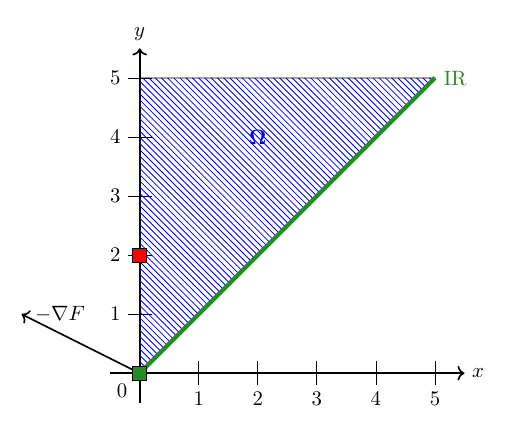}
    \caption{An illustration of the bilevel problem~\eqref{eq:ExampleBilevelProb_Optimistic}: the set $\Omega$ (blue), the inducible region IR (green line), the optimal solution $(0,0)$ (green square), and the non-admissible point $(0,2)$ (red square).}
    \label{fig:ExampleNeedAdmissibility}
\end{figure}

Now assume that a solver, applied to the lower-level problem at $x^0 = 0$, erroneously returns $y^0 = 2$ as the optimal lower-level response.
As shown by the red square in~\Cref{fig:ExampleNeedAdmissibility}, the point $(x^0,y^0) = (0,2)$ belongs to $\Omega$ but not to IR: it satisfies the lower-level constraints, but $y^0$ is not lower-level optimal.
Nevertheless, its upper-level objective function value $F(0,2) = -2$ is strictly less than the true optimal value $F(0,0) = 0$.
A BL-DFO algorithm reporting this point would appear to outperform an algorithm that correctly returns $(0,0)$, although its solution is not admissible.
This example shows how an inaccurately solved lower-level problem can make an algorithm appear better than it actually is.
It also highlights the importance of verifying the optimality of the lower-level response associated with any upper-level variable~$x$.}

When the lower-level response is generated {\bl inexactly}, it is possible that it is not optimal for the lower-level due to local minima or an inability to converge.
Therefore, that solution should have been discarded by the upper-level decision maker since non-admissible solutions cannot be considered optimal.
Consequently, it is possible that an inefficient lower-level algorithm produces falsely admissible solutions. 
These solutions could satisfy the upper-level constraints and have an excellent upper-level objective function value.
To summarize, this inefficient algorithm would claim to produce excellent solutions to the bilevel problem.

In this context, a procedure that challenges the admissibility of generated points, named the \textit{refereeing procedure}, is required to ensure a fair comparison of BLO algorithms in a DFO setting.
The goal of such a procedure is to detect non-admissible points and remove them from the benchmarking data. 
Points where the lower-level responses are computed by an inefficient algorithm are therefore likely to be discarded.
A referee takes as input the run log of an optimization run, analyzes entries that are claimed to be admissible, and flags the ones that are shown to be non-admissible.
The revised log is then used by the benchmarking process.
{\bl Ensuring lower-level optimality of bilevel solutions goes beyond the framework of DFO; refereeing procedures could be applied to any context where computing the global lower-level solution for each $x$ is impossible.}

Additionally, solving the lower-level problem for a given upper-level variable~$x$ may be computationally expensive.
{\bl Although CPU time could be used to measure this effort, DFO benchmarking techniques~\citep{DoMo02,MoWi2009} rely on the number of function evaluations, the main budget unit of derivative-free algorithms.
Accordingly, we measure the total computational effort of BL-DFO algorithms by the number of upper- and lower-level function evaluations, although some of the tools proposed here extend to CPU time.}

The computational effort devoted to {\bl solve the lower-level problem} by BL-DFO algorithms is not always considered by existing benchmarking techniques, while it represents a significant expense of derivative-free methods in solving a bilevel problem.
Besides, applications considered may have upper- and lower-level problems whose evaluation costs differ.
A scaled computational effort that aggregates both upper- and lower-level evaluations allows then to account for the general cost involved in solving~\eqref{eq:MainBilevelProb_Optimistic}.
This computational effort can then be incorporated into existing benchmarking techniques to measure the true expense employed to solve a bilevel problem.

The contributions of this work are summarized as follows.

\begin{itemize}
    \item Notions about feasibility and admissibility in BL-DFO are extended to an inexact setting to cover practical issues.
    In addition, a review of existing BL-DFO algorithms is conducted, focusing on the inexactness of the methods employed to solve the lower-level problem and the benchmarking techniques used to compare such algorithms.
    \item A refereeing procedure is introduced to challenge the admissibility of points through additional optimization runs.
    Several strategies and settings are suggested to devise multiple refereeing strategies.
    \item A scaled computational effort that aggregates upper- and lower-level evaluations with a scaling factor is introduced.
    The combination of the scaled computational effort and the refereeing procedure yields an extension of single-level benchmarking techniques that accounts for the specific issues inherent to BLO.
    \item The proposed methodology is illustrated on a collection of analytical problems from the BOLIB library~\citep{Zhou2021}, adapted to the Julia programming language, and available as an open source package entitled \texttt{BOLIB.jl}\footnote{\href{https://github.com/Vanadjy/BOLIB.jl}{https://github.com/Vanadjy/BOLIB.jl}}.
    The refereeing strategies are analyzed along with tests about the scaled computational effort to assess their impact on data profiles~\citep{MoWi2009}.
\end{itemize}

The paper is structured as follows. 
\Cref{sec:BLO_Inexact_DFO} explains how BLO can be considered in inexact and DFO settings.
BLO notions are extended to the inexact setting, followed by an extensive literature review of BL-DFO algorithms that employ both exact and inexact methods to solve the lower-level problem.
The refereeing procedure is detailed in~\Cref{sec:CompareBLDFOalgs} along with three refereeing strategies.
Possible settings for the refereeing procedure are also discussed.
The scaled computational effort is introduced in the description of the methodology as well.
These features result in an extension of existing benchmarking techniques to BLO.
Their impact is illustrated on numerical examples in~\Cref{subsec:ComputationalTests}.

\paragraph{Notations:} $\R_+ \coloneqq [0,+\infty)$ denotes the set of nonnegative real numbers, and $\overline{\R} \coloneqq \R \cup \{\pm \infty\}$ represents the extended real numbers. 
For any integer \(n \ge 1\), \(\intInterval{1}{n}\) denotes the set \(\{1,2,\dots,n\}\).

\section{Inexact bilevel optimization setting and DFO}\label{sec:BLO_Inexact_DFO}

BL-DFO methods studied in this work proceed iteratively by proposing an upper-level variable $x$, and then computing the associated lower-level response $y$.
However, algorithms employed to solve the lower-level problem may produce inexact lower-level responses with respect to a tolerance $\varepsilon$.
Setting $\varepsilon = 0$ means the lower-level problem is solved exactly.
Additionally, bilevel optimization requires a dedicated terminology to characterize points that satisfy lower-level optimality and/or inequality constraints.
This bilevel specific terminology is then extended to an inexact framework. 
A review of BL-DFO algorithms is presented to identify if the methods employed at the lower-level solve it exactly. 
Commonly used benchmarking techniques are also identified in this literature review and are further discussed. 

\subsection{BLO in an inexact setting}\label{subsec:InexactBLO}

{\bl The specific structure of bilevel problems requires a specific terminology~\citep{HaJaSa91} to characterize solutions produced by BL-DFO algorithms.}
This terminology needs to be extended to an inexact setting for situations in which lower-level responses associated to an upper-level variable $x$ are computed inexactly.

A point $(x,y) \in X \times Y$ is said to be \textbf{feasible} for the bilevel problem~\eqref{eq:MainBilevelProb_Optimistic} if {\bl $y \in Y(x)$ and $G(x,y) \le 0$}.
The set of all feasible points is denoted by $\Omega$.
The set of optimal lower-level {\bf responses} associated with an upper-level variable $x \in X$ is denoted by
\begin{equation}\label{eq:PsiDef}
    \Psi(x) \coloneqq \argmin {y \in Y(x)} f(x,y).
\end{equation} 
A point $(x,y) \in X \times Y$ is said to be \textbf{rational} if $y$ is an optimal lower-level response associated with $x$, i.e., $y \in \Psi(x)$.
A point \((x,y) \in X \times Y\) is said to be \textbf{admissible} if \((x,y) \in \Omega\) and \(y \in \Psi(x)\), i.e., $(x,y)$ is feasible and rational.
The set of all admissible points is called the \textbf{induced region}, denoted by \(\IR\).
Finally, a point $(x,y) \in X\times Y$ is \textbf{optimal} for~\eqref{eq:MainBilevelProb_Optimistic} if $(x,y) \in \IR$ and
\[
F(x,y) \le F(\bar x, \bar y), \quad \forall (\bar x, \bar y) \in \IR.
\]
BLO problems are inherently difficult.
Even if all functions are linear,
 the induced region can be nonconvex and disconnected,
 and finding a global optimal solution {\bl can be strongly NP-hard}~\cite{HaJaSa91}. 

Some algorithms typically compute approximate minimizers of the lower-level optimization problem. 
To accommodate this inexactness, the notions of inexact {\bl rationality and admissibility} are formalized in~\Cref{def:LL_eps_optimal}.

{\bl \begin{definition}\label{def:LL_eps_optimal}
    Let $\varepsilon \in \R_+$ be a tolerance threshold. 
        Consider an upper-level variable $x \in X$ and its associated lower-level response $y_{\varepsilon} \in Y(x)$.
        Then, the point $(x, y_{\varepsilon})$ is said to be $\varepsilon$\textbf{-rational} if
        \begin{equation*}\label{eq:ConditionEpsRational}
            f(x,y_{\varepsilon}) \le  f(x,y) + \varepsilon, \quad \textrm{ for all } y \in Y(x).
        \end{equation*}
        The point $(x, y_{\varepsilon})$ is said to be $\varepsilon$\textbf{-admissible} if it is feasible and $\varepsilon${-rational}.
\end{definition}}


BL-DFO algorithms that employ inexact methods to compute lower-level responses are commonly used.
Extending BLO terminology to the inexact setting therefore enables consideration of a relaxed criterion to revoke the admissibility of generated solutions. In this case, inexact methods used to solve the lower-level problem are not unfairly penalized in the benchmarking process.

{\bl \subsection{A review of BL-DFO algorithms and benchmarking techniques}\label{subsec:BLO_Inexact_DFO_Algs}}
Various benchmarking techniques are available to compare BL-DFO methods, but  are often applied unevenly in the literature.
This section provides an extensive literature review of existing BL-DFO methods to highlight the benchmarking techniques used in bilevel optimization.
In addition to the benchmarking strategies,
the review lists methods that solve the lower-level problem and mentions if it is assumed to be exact.
All algorithms are classified based on the DFO scheme applied to the upper-level problem.

\paragraph{Direct search methods:}
A common approach for solving blackbox optimization problems is to exploit information obtained directly from objective function evaluations; hence the name direct-search {\bl (DS)} methods~\cite{HoJe61a}.
These approaches either employ a mesh to select points and consider a simple decrease acceptance condition, or they consider a sufficient decrease criterion to accept new iterates~\cite{LeToTr00a}. 
More recently, \citet{DzRiRoZe2025}
 propose a unified review of the convergence analyses of direct search methods.

In the bilevel setting, \citet{MeDe2011} developed a direct search framework along with a sufficient decrease condition.
Two direct search variants for bilevel programming are proposed.
The first requires that the set of poll directions has a sufficiently large vector density, while the second changes the set of poll directions for each unsuccessful iteration.
The lower-level response is {\em a priori} given by an exact oracle for each upper-level variable $x$. \citet{MeDe2011} proposed and analyzed directional direct-search methods in the BLO context; however, no implementation details for the proposed algorithms were provided.
The first adaptation of a mesh adaptive direct-search~\cite{AuDe2006} scheme was analyzed by~\citet{DiKuRiZe2024}, where each lower-level response is obtained from an inexact oracle. They also studied direct-search strategies based on sufficient decrease in both smooth and nonsmooth settings.  The performance of their methods is assessed using data profiles based on upper-level function evaluations.

For health insurance optimization problems,~\citet{ZhLi2014} developed direct search schemes to solve both upper-level and lower-level problems.
The lower-level responses are obtained by applying a direct search scheme with sufficient decrease, and the tolerance for stopping the lower-level optimization is updated at each iteration.
For one of their applications, the lower-level problem is solved using a variant of a trust-region method tailored to Nash-equilibrium problems~\citep{Yuan2011}.
They illustrate this approach with convergence plots on problems of hospital competitions~\citep{EgYi2004} and oligopolistic market models.

\paragraph{Line-search methods:} 
Recently, a BL-DFO framework was proposed by~\citet{Cesaroni2026}, which includes a line search {\bl (LS)} strategy performed by a projected extrapolation procedure that computes a direction $d_k$ and a step length $\alpha_k$.
The upper-level variable is updated to decrease a hyper-objective function whose lower-level responses are computed with an inexact oracle.
The performance of the proposed method is illustrated using tables, performance profiles, and data profiles.

\paragraph{Trust-region methods:}
Surrogates of the objective function using gradient and Hessian estimates can be constructed to apply trust-region {\bl (TR)} methods in a derivative-free context.
These approaches are adapted to problems where the derivatives exist but are not accessible.
\citet{CoMaSa2005} adapted trust-region schemes to bilevel optimization, and~\citet{CoVi2012} extended this framework to the DFO setting.
They minimize quadratic interpolation models (minimum Frobenius norm models when not enough previously computed points are available) for both upper-level and lower-level problems.
The lower-level problem is allowed to be solved inexactly, and the construction of models may exploit previously evaluated points.
Their tests are conducted on analytical problems with linear box constraints for some of them.
Convergence plots are drawn, showing the upper-level function value with respect to upper-level and lower-level function evaluations separately.

The algorithm proposed by~\citet{EhRo2021} uses a trust-region scheme to solve the upper-level problem, while lower-level responses are obtained via an inexact oracle.
They use either a gradient descent method or the fast iterative shrinkage-thresholding algorithm (FISTA)~\citep{BeTe2009a} to solve a strongly convex lower-level problem arising from image denoising applications.
The upper-level problem consists of a regularized non-linear least squares structure addressed in a DFO setting.
This method is illustrated on mathematical imaging problems through convergence plots showing the upper-level value function with respect to lower-level evaluations.

\paragraph{Finite-differences methods:}
BL-DFO algorithms may exploit estimates of gradients or Hessians {\bl through finite differences (FD)} to drive the optimization process.
\citet{AgGh2025} estimate an upper-level Hessian inverse and partial gradient approximations using i.i.d.~random vectors drawn from a standard Gaussian distribution.
The lower-level problem is solved by a derivative-free stochastic gradient descent method, yielding an inexact lower-level response, which is then used in the outer optimization loop.
The upper-level variables are also updated through a derivative-free projected stochastic gradient method.
Computational tests are conducted on parameter optimization for classification models and display convergence plots with respect to both upper- and lower-level function evaluations separately.

{\bl \citet{Staudigl2025} propose an alternative stochastic derivative-free proximal-gradient framework that accommodates both exact and inexact lower-level responses.
The framework is tested on signal denoising and optimal experimental design problems.
 Performance is assessed through plots of the normalized empirical error of the lower-level solutions, together with application-specific figures, rather than through a common benchmarking methodology.} 
For transportation network problems,  \citet{Maheshwari2023} employ a derivative-free gradient descent to update the upper-level variables and a projected gradient descent to compute inexact lower-level responses. Numerical experiments are {\bl illustrated} by convergence plots.

\paragraph{Bayesian optimization methods:}
Another approach to solve problems in DFO is based on the training of a Gaussian process (GP) on a design of experiments (DoE).
This approach, known as Bayesian optimization~(BO), consists in optimizing an acquisition function that yields trial points to enrich the DoE used to update the GP.
Different strategies are employed for either exploration or intensification purposes with respect to the selected acquisition function and the values of its associated parameters.

To adapt BO to bilevel optimization,~\citet{Ekmekcioglu2024} design the algorithm BILBAO, which constructs a GP for both upper- and lower-level problems.
The Regional Expected Value of Improvement constitutes the lower-level acquisition function and is optimized using a quasi-Newton method combined with a multi-start strategy.
The resulting point $(x,y)$ is used to update the lower-level GP, and a mapping function $\Phi$ is used to update upper-level variables, incurring inexact lower-level responses.
Convergence plots related to the inexactness of lower-level responses are reported with respect to upper-level and lower-level function calls combined, assuming their evaluation costs are identical.
Authors in~\cite{Kieffer2017,DoPr2023} present similar methods in which the upper-level variables are computed through the optimization of the lower confidence bound acquisition function.
Lower-level responses are computed by a sequential least squares programming~(SLSQP) method~\citep{Kraft1988}.
{\bl Results are reported in tables of acquisition function fitness scores and average numbers of function calls, complemented by convergence and trajectory plots of the acquisition functions.} 

To combine the advantages of several acquisition functions,~\citet{DoPr2023-2} solve a multi-objective problem to suggest upper-level trial variables.
The lower-level structure consists of a least squares problem and is solved using the SLSQP method.
{\bl The convergence plot of the difference between the true optimal upper-level value and the approximated one is displayed.}
The effort employed concerns upper-level evaluations (outer iterations) only.

Recently, the BILBO approach devised by~\citet{Chew2025} optimizes the upper confidence bound~(UCB) acquisition function to solve the upper-level problem.
The lower-level function $f$ is minimized through the UCB acquisition function only if $f$ constitutes the largest regret bound.
The latter is computed from an estimated lower-level optimal response on a trusted set of feasible solutions.
The lower-level responses are supposedly computed exactly, but this can be extended to an inexact setting.
Results are reported using convergence plots with respect to evaluations of any function involved in the bilevel problem, showing the decrease of the instant regret, which represents the loss associated with not selecting the optimal point.
Based on an application of blackbox attacks on finite sequences,~\citet{Man2025} build a multi-kernel GP whose weights are updated by solving the upper-level problem block-wise.
The acquisition function employed to solve the lower-level problem is the Expected Improvement, and the GP parameters are updated using maximum {\em a posteriori} estimations.

\paragraph{Heuristic-based methods:}
{\bl Heuristic DFO methods, despite lacking convergence guarantees, can perform well on certain problems.}
These methods can be used as standalone solvers or to enhance the performance of existing algorithms in a DFO setting.
Numerous existing BL-DFO algorithms employ {\bl evolutionary algorithms (EA)} to solve the upper-level or the lower-level problems, notably genetic algorithms (GA).

The method developed by~\citet{Jiang2023} solves both upper- and lower-level problems by mutating populations of solutions using a simulated binary crossover and a polynomial mutation.
A part of the upper- and lower-level populations is selected by a surrogate Kriging model, which approximates the relationship between upper-level variables and their associated lower-level responses.
Lower-level optimization is performed by a GA with dynamic accuracy for each upper-level variable.
Tables representing upper-level and lower-level accuracy, evaluation times, and Euclidean distance to the best solution  are used to assess the performance of their method.

The approach proposed by~\citet{MaHu2025} addresses the upper-level problem with a Differential Evolution~(DE)~\citep{DaSu2010} strategy.
The lower-level problem represents a surrogate learning stage performed by two successive gradient descent methods on two different loss functions, yielding a  Kolmogorov-Arnold Network~\citep{Liu2024} used to approximate the lower-level objective function.
Tables of mean scores and cost convergence plots with respect to lower-level evaluations are displayed to compare the different existing methods and their variants.
The method devised by~\citet{Islam2017} uses a DE algorithm for the upper-level problem, while the lower-level function is optimized by another DE through three different surrogate models: linear and quadratic regressions, and Kriging models.
Upper-level function values are displayed through tables and convergence plots for two different sets of instances~\citep{Sinha2014} with respect to the total number of function evaluations (upper-level and lower-level combined with equivalent evaluation costs).
Performance profiles are constructed from this set of instances, with an accuracy value defined by median accuracy and total function evaluations.
The method for calculating the accuracy value is not mathematically explicit, and the admissibility of points is not discussed.

In a similar work of~\citet{Islam2018}, a DE method is used for the upper-level framework, where each lower-level response is obtained from a BO algorithm~\citep{JoScWe1998}.
An interior point algorithm is applied after BO to ensure the optimality of lower-level responses.
The results of their tests are summarized in tables.
Convergence plots showing the evolution of the upper-level function value with respect to upper-level and lower-level evaluations, equally combined, are also displayed.
 For a human gait application problem, \citet{Nguyen2019} propose a GA framework for the upper-level problem, while a collocation strategy is used to compute lower-level responses by calling the IPOPT solver~\citep{ipopt}.
Their tests are illustrated by both tables showing solution errors and convergence plots over the number of generations of the GA.

{\bl \Cref{table:ArticlesSummarizeBLDFO} summarizes this literature review, specifying for each reference the upper-level methodology, the lower-level solver and whether it is exact, and the benchmarking tools employed.
The first column contains the source reference; the second indicates the upper-level algorithmic strategy.
Then a checkmark indicates that the lower-level exact solutions can be obtained.
The wide central column indicates the lower-level algorithmic strategy.
Finally, the last three columns indicate whether the proposed analysis relies on tables, convergence plots, and/or profiles.}

\begin{table}[htb!]
\centering
{\footnotesize \renewcommand{\tabcolsep}{2.1pt}
\begin{tabular}{|l|l|cl|l|l|l|}
    \hline
   Reference                             & UL & {\bl Claimed Exact} & LL Solver & T & C & P \\ \hline
    \citet{MeDe2011}     &   DS     &  \OK  & Unspecified     &       &         &   \\
    \citet{DiKuRiZe2024}     &   DS     &   &  Inexact Oracle
    &       &         &  \OK \\  
    \citet{ZhLi2014}     &   DS  &   &   
DFO method (direct search or trust-region)   
    &   \OK    &    \OK     &   \\
    \citet{Cesaroni2026}     &   LS  &   &   
Inexact Oracle
    &   \OK    &        &  \OK  \\
\citet{CoVi2012}               &   TR    &     &     Trust-region method       &       &   \OK   &   \\
\citet{EhRo2021}      &   TR    &    &     Gradient descent method or FISTA 
  &       &   \OK   &   \\
    \citet{AgGh2025}    &   FD    &      &  Gradient descent method &      &   \OK   &  \\
    \citet{Staudigl2025}    &   FD   & \OK   &   Exact and inexact Oracle
    &      &  \OK   &  \\
    \citet{Maheshwari2023}    &   FD  &    &    Projected gradient method        &      &   \OK  &  \\
\citet{Ekmekcioglu2024}    &   BO    &     &   Acquisition function-based BO method    &       &   \OK   &  \\
    \citet{DoPr2023}    &   BO      &   \OK & Nonlinear problem solved by SLSQP            &    \OK   &     &  \\
    \citet{Kieffer2017}    &   BO      &  \OK & Nonlinear problem solved by SLSQP  &    \OK   &     &  \\
    \citet{DoPr2023-2}    &   BO      &   \OK & Least squares problem solved by SLSQP  &  \OK  &  \OK   &  \\
         \citet{Chew2025}    &   BO  &     &  Acquisition function-based BO method &    &  \OK   &  \\ 
    \citet{Man2025}    &   BO   &   &  Acquisition-function based BO method  &  \OK  &  \OK   &  \\ 
    \citet{Jiang2023}    &   EA   &   &  Population-based algorithm   &  \OK  &    &  \\
    \citet{MaHu2025}    &   EA      &   & 
    Gradient descent method  &  \OK  &  \OK  &  \\
    \citet{Islam2017}    &   EA   &    &    Surrogate-assisted  differential evolution  method &  \OK  &  \OK  & \OK \\
    \citet{Islam2018}    &   EA   & \OK    &   Two phases: a BO method then an interior point method  &  \OK  &  \OK  &  \\
    \citet{Nguyen2019}    &   EA &  \OK     & Direct collocation method     &   \OK   &   \OK  &  \\ \hline
\end{tabular}
}
\caption{Existing BL-DFO algorithms, with emphasis on upper-level (UL) and lower-level (LL) methodologies. The last three columns indicate if the analysis uses tables (T), convergence plots (C), or benchmarking profiles (P).}
\label{table:ArticlesSummarizeBLDFO}
\end{table}

\subsection{Limits of existing benchmarking strategies in bilevel setting}\label{subsec:ExistingStrats}

The references above suggest extensions of single-level benchmarking techniques to a bilevel framework.
To this end, the $\varepsilon$-admissibility of points generated by BL-DFO algorithms is assumed to be guaranteed.
Yet, inexactly computed lower-level responses may result in non-admissible (and {\em a fortiori} non-optimal) solutions.
Even claimed exact methods that solve the lower-level problem may not find a global optimal solution, as this would require expensive global optimization tools.
Challenging the $\varepsilon$-admissibility of generated solutions is therefore required to extend benchmarking techniques to the bilevel setting to ensure fair comparisons between BL-DFO algorithms.

Existing benchmarking techniques account for either upper-level evaluations~\citep{DiKuRiZe2024,DoPr2023-2,CoVi2012} or lower-level evaluations~\citep{Cesaroni2026,EhRo2021,Ma2025,CoVi2012} as cost metrics.
Counting upper-level evaluations alone disregards the computational effort involved in solving the lower-level problem, which constitutes a substantial cost in solving~\eqref{eq:MainBilevelProb_Optimistic}.
Conversely, restricting the cost metric to lower-level evaluations is inadequate when evaluating the upper-level involves a significant expense.
Both upper-level and lower-level evaluations are aggregated in some works in the literature (see~\citep{Ekmekcioglu2024,Islam2017,Islam2018}), assuming that upper- and lower-level evaluation costs are equivalent.
Introducing a scaling factor that reflects the relative cost between the two levels is required to properly assess the total computational effort deployed to solve the bilevel problem.

Benchmarking techniques used to compare existing BL-DFO algorithms include tables, convergence plots, and benchmarking profiles, none of which fully address the difficulties arising from BLO.

\begin{itemize}
    \item Tables reporting solving statistics become difficult to interpret at scale, making them unsuitable for systematic comparison across large algorithm and instance sets.
    \item Convergence plots illustrate the evolution of algorithms on individual instances and, thus, hardly support conclusions drawn over a broad collection of problem instances.
    \item {\bl Benchmarking profiles compare algorithms by plotting the objective function values of feasible points against the computational effort required to obtain them.
In BLO, however, benchmarking must also account for the $\varepsilon$-admissibility of these points and for the computational effort at both the upper- and lower-levels.
    }
\end{itemize}

This review of BL-DFO methods motivates the need to generalize the methodology of current benchmarking techniques to the bilevel setting.
To this end, a refereeing procedure is introduced to describe how to challenge the $\varepsilon$-admissibility of points, and a scaling method is proposed to properly aggregate upper-level and lower-level evaluations into a unified effort metric.

\section{A referee-based bilevel benchmarking technique}\label{sec:CompareBLDFOalgs}

Current benchmarking techniques do not adequately consider the following BLO issues.
Firstly, the $\varepsilon$-admissibility of solutions produced by a BL-DFO algorithm is not always guaranteed. 
Challenging the $\varepsilon$-admissibility of points through a procedure that discards non-admissible points is therefore required.
Secondly, upper- or lower-level evaluations are either ignored or not properly aggregated.
A scaled computational effort that accounts for the actual expenses at both levels is therefore needed to reflect the general cost to solve problem~\eqref{eq:MainBilevelProb_Optimistic}.
A refereeing procedure is introduced along with refereeing strategies to challenge the $\varepsilon$-admissibility of solutions. 
Possible settings for the refereeing procedure are then detailed.
Finally, a scaled computational effort is introduced and the evaluation budget units used in data profiles are discussed to extend this benchmarking technique to BLO.

\subsection{A refereeing procedure to revoke admissibility}
\label{subsec:RefAdmissibility}

{\bl In order to challenge the claimed $\varepsilon$-admissibility of points produced by a BL-DFO algorithm, a post-processing \emph{refereeing procedure} is introduced.} 
A criterion of non-admissibility is first given to revoke solutions that are not $\varepsilon$-admissible before detailing the procedure.
For a given upper-level variable \(x \in X\), the corresponding lower-level response $y \in Y$ may be suboptimal, as the lower-level algorithm might terminate prematurely, converge to a local minimizer, or simply fail.
In blackbox, nonsmooth or nonconvex settings, global optimality is generally unattainable.
The $\varepsilon$-admissibility of a feasible point \((x,y) \in X \times Y\) is revoked when an alternative feasible point \((x, \hat y)\in X \times Y\) is shown to satisfy
\[f(x,\hat y) < f(x, y) - \varepsilon.\]

Consider a benchmark set of instances \(\Pbank\) and a collection of algorithms \(\A\).
The instances of \(\Pbank\) should be chosen to be as general as possible by avoiding biased initial conditions or artificial structures.
Additionally, it should be representative of the class of target problems of interest.
Let $\mathcal{H}_{a,p}$ be the set of claimed $\varepsilon$-admissible points generated by applying algorithm $a \in \A$ on the problem instance $p \in \Pbank$.
Let $H \in \N$ be the cardinality of $\mathcal{H}_{a,p}$, so that this set can be written as 
\begin{equation}\label{eq:Hap}
    \mathcal{H}_{a,p} \coloneqq \{(x_k,y_k) \,:\, k \in \intInterval{0}{H}\}.
\end{equation}
\Cref{def:referee} describes the {\em a posteriori} procedure to challenge, and eventually to revoke, the $\varepsilon$-admissibility of a point of $\mathcal{H}_{a,p}$.

\begin{definition}\label{def:referee}
    Let $\varepsilon \in \R_{+}$ be a tolerance threshold and 
    $(x, y) \in \mathcal{H}_{a,p}$,
    a feasible point for problem instance $p$ produced by the BLO algorithm $a$.
    A {\bf referee} is an algorithm $r$ that {\bf challenges} the $\varepsilon$-admissibility of $(x, y)$ through the following steps:
    \begin{enumerate} 
    \item Apply the referee $r$ with  upper-level variable $x$ to produce a lower-level response $y^r$ ;
    \item If $(x,y^r)$ is feasible and if $f(x,y^r) < f(x, y) - \varepsilon$,
     then the $\varepsilon$-admissibility of $(x, y)$ is {\bf revoked}. 
    \end{enumerate}
\end{definition}

Each referee can belong to $\A$ (Internal Referee) or not (External Referee).
A referee does not certify the admissibility of a point \((x,y)\).
Instead, it challenges possibly non-admissible points by looking for an alternative lower-level response that contradicts the presumed $\varepsilon$-lower-level optimality of the response \(y\) for the same \(x\).
The refereeing procedure is intended to be applied as a preprocessing step to the benchmarking process.

Let $r$ be a referee.
The refereeing procedure may be specified differently depending on the purposes of the optimization framework.
Accordingly, three distinct refereeing strategies are defined, each tailored to a specific validation criterion, yielding a filtered history $\mathcal{H}_{a,p}^{r} \subseteq \mathcal{H}_{a,p}$.
\begin{itemize}
    \item The End-Point Referee challenges only the $\varepsilon$-admissibility of the final point $(x_H, y_H)$.
    \item The Complete Referee challenges the $\varepsilon$-admissibility of all points $(x_k, y_k)$ for $k \in \intInterval{0}{H}$.
    \item The Reverse Referee challenges the $\varepsilon$-admissibility of some points in reverse order, starting from the final point $(x_H, y_H)$.
\end{itemize}

The End-Point Referee challenges only the last point of \(\mathcal{H}_{a,p}\), namely $(x_H,y_H)$.
If the $\varepsilon$-admissibility of $(x_H,y_H)$ is revoked, then the entire execution history $\mathcal{H}_{a,p}$ is discarded.
Otherwise, all the points \((x_k, y_k)\) for $k \in \intInterval{0}{H}$ are preserved in $\mathcal{H}_{a,p}^r$.
Then, for the End-Point strategy,
\[
\mathcal{H}_{a,p}^r = \begin{cases}
    \mathcal{H}_{a,p} & \textrm{if } (x_H,y_H) \textrm{ is not revoked by } r, \\
    \emptyset & \textrm{otherwise.}
\end{cases}
\]
This strategy disregards the trajectory of the algorithms and bases its decision solely on the $\varepsilon$-admissibility of the returned solution. 
The End-Point Referee is inherently myopic, as it may discard points that might be $\varepsilon$-admissible simply because the last point is non-admissible.

The Complete Referee challenges each point $(x_k, y_k) \in \mathcal{H}_{a,p}$ and discards every revoked point from \(\mathcal{H}_{a,p}\) yielding the corresponding filtered history $\mathcal{H}_{a,p}^{r}$, which represents $\mathcal{H}_{a,p}$ deprived of the points revoked by $r$, i.e.,
\[
\mathcal{H}_{a,p}^r = \left\{ (x_k,y_k) \in \mathcal{H}_{a,p} \,:\, 
    k \in \intInterval{0}{H} \mbox{ and }(x_k,y_k) \mbox{ is not revoked by } r\right\}.
\]
The resulting filtered history $\mathcal{H}_{a,p}^{r}$ can then be used to construct performance and data profiles~\citep{DoMo02,MoWi2009}.
However, the Complete Referee is computationally expensive, as it requires rerunning multiple optimization procedures for each history $\mathcal{H}_{a,p}$.

The Reverse Referee is designed to mitigate the drawbacks of the End-Point and the Complete Referees.
It challenges the $\varepsilon$-admissibility of the last point \((x_H,y_H)\) first, then proceeds backward through $\mathcal{H}_{a,p}$ until either an $\varepsilon$-admissible point is identified, denoted by $(x_{K},y_{K})$, or the $\varepsilon$-admissibility of the  point \((x_0,y_0)\) is revoked.
If such an integer $K$ exists, then all points $(x_k,y_k)$ for $k \in \intInterval{0}{K}$ are preserved, thus
\[
\mathcal{H}_{a,p}^r =  \left\{   \renewcommand{\arraycolsep}{2pt}
     \begin{array}{rl} 
     (x_k,y_k) \in \mathcal{H}_{a,p} \,:\, 
     & k  \in \intInterval{0}{K}  \mbox{ where } K \in \intInterval{0}{H}  \mbox{ satisfies} \\
     & (x_{K},y_{K})  \mbox{ is not revoked by } r,\\
     & (x_{\ell},y_{\ell})  \mbox{ is revoked by } r 
            \mbox{ for all } \ell \in \intInterval{K +1}{H}
\end{array} \right \}.\\
\]
Otherwise, all points of $\mathcal{H}_{a,p}$ are revoked, resulting in $\mathcal{H}_{a,p}^r = \emptyset$.

If the refereeing procedure is deterministic, then both the Complete and Reverse referees identify the same final $\varepsilon$-admissible solution.
The Reverse Referee is therefore well-suited for constructing accuracy profiles~\citep{Beiranvand2017}, as it provides the final $\varepsilon$-admissible solution for each algorithm-instance point and avoids verifying the entire history.

{\bl To illustrate the outcomes of each refereeing strategy, consider two hypothetical BL-DFO algorithms, \texttt{A} and \texttt{B}, each producing a history of points.
For each algorithm, \Cref{table:IllustrationRefereeOutcomes} reports the number of upper-level function evaluations together with the best objective function value claimed so far.
\Cref{table:data} displays the original data, with revoked points in red, and \Cref{table:EndPointReferee,table:ReverseReferee,table:CompleteReferee} show the data filtered by each of the three refereeing strategies.
This small example shows that the choice of referee can change the outcome of a comparison.
With the Reverse referee, algorithm~\texttt{B} requires $123$ evaluations to $2\%$-solve the problem, whereas it requires $200$ evaluations with the Complete referee.
With the End-Point referee, algorithm~\texttt{A} fails to $2\%$-solve the problem.
}

\begin{figure}[htb!]
    \centering
    \footnotesize
    \setlength{\tabcolsep}{4pt}
    \subcaptionbox{Original data\label{table:data}}
    {\begin{tabular}[t]{|c|c||c|c|}
        \hline
        $N_A$ & $F_A$ & $N_B$ & $F_B$ \\ \hline
        $1$ & $695$ & $1$ & $695$ \\
        {\color{red} $42$} & {\color{red} $562$} & $5$ & $646$ \\
        $55$ & $506$ & {\rd $10$} & {\rd $618$} \\
        $68$ & $212$ & {\rd $20$} & {\rd $555$} \\
        $104$ & $100$ & $80$ & $499$ \\
        {\rd $152$} & {\rd $58$} & {\rd $123$} & {\rd $102$} \\
        {\rd $200$} & {\rd $9$} & $200$ & $101$ \\ \hline
    \end{tabular}}
    \hfill
    \subcaptionbox{End-Point Referee\label{table:EndPointReferee}}
    {\begin{tabular}[t]{|c|c||c|c|}
        \hline
        $N_A$ & $F_A$ & $N_B$ & $F_B$ \\ \hline
        $1$ & $695$ & $1$ & $695$ \\
        - & - & $5$ & $646$ \\
        - & - & $10$ & $618$ \\
        - & - & $20$ & $555$ \\
        - & - & $80$ & $499$ \\
        - & - & $123$ & $102$ \\
        - & - & $200$ & $101$ \\ \hline
    \end{tabular}}
    \hfill
    \subcaptionbox{Reverse Referee\label{table:ReverseReferee}}
    {\begin{tabular}[t]{|c|c||c|c|}
        \hline
        $N_A$ & $F_A$ & $N_B$ & $F_B$ \\ \hline
        $1$ & $695$ & $1$ & $695$ \\
        $42$ & $562$ & $5$ & $646$ \\
        $55$ & $506$ & $10$ & $618$ \\
        $68$ & $212$ & $20$ & $555$ \\
        $104$ & $100$ & $80$ & $499$ \\
        - & - & $123$ & $102$ \\
        - & - & $200$ & $101$ \\ \hline
    \end{tabular}}
    \hfill
    \subcaptionbox{Complete Referee\label{table:CompleteReferee}}
    {\begin{tabular}[t]{|c|c||c|c|}
        \hline
        $N_A$ & $F_A$ & $N_B$ & $F_B$ \\ \hline
        $1$ & $695$ & $1$ & $695$ \\
        - & - & $5$ & $646$ \\
        $55$ & $506$ & - & - \\
        $68$ & $212$ & - & - \\
        $104$ & $100$ & $80$ & $499$ \\
        - & - & - & - \\
        - & - & $200$ & $101$ \\ \hline
    \end{tabular}}
    \caption{An illustration of the outcomes of the refereeing strategies on two BL-DFO algorithms.}
    \label{table:IllustrationRefereeOutcomes}
\end{figure}

{\bl The three refereeing strategies serve distinct purposes.
The End-Point referee verifies the $\varepsilon$-admissibility of the last point computed by a BL-DFO algorithm.
It can serve as a preliminary benchmark, indicating which methods return a solution that is not $\varepsilon$-admissible.
The Reverse referee determines the point beyond which the remainder of the trajectory is not $\varepsilon$-admissible.
Finally, the Complete referee verifies the $\varepsilon$-admissibility of every point generated by the algorithm, thereby identifying which parts of its trajectory are not $\varepsilon$-admissible.
}

Several algorithms may be used for a single refereeing procedure. 
Denote by $\Refs$ the set of referees selected to challenge the $\varepsilon$-admissibility of points of $\mathcal{H}_{a,p}$.
The filtered history resulting from this refereeing procedure is denoted by
\begin{equation}\label{eq:Hap_filtered}
    \mathcal{H}_{a,p}^{\Refs} \coloneqq \bigcap_{r \in \Refs} \mathcal{H}_{a,p}^r.
\end{equation}

Refereeing strategies differ in computational cost as they do not challenge the same number of points of $\mathcal{H}_{a,p}$.
The End-Point Referee is the least expensive strategy since it challenges only the last point generated during the run of algorithm $a$ on instance $p$.
Its overall cost is therefore proportional to the number of instances $|\Pbank|$.
The Complete Referee challenges the $\varepsilon$-admissibility of all points in $\mathcal{H}_{a,p}$, so the cost of this refereeing process scales with the length of each execution history.
The Reverse Referee has a lower computational cost than the Complete Referee, as it challenges the $\varepsilon$-admissibility on a subset of $\mathcal{H}_{a,p}$, proceeding backward from the final point $(x_H, y_H)$.
In the worst case scenario, i.e., when the $\varepsilon$-admissibility of each point is revoked, its computational cost coincides with that of the Complete Referee.
Conversely, the Reverse Referee cost is identical to that of the End-Point Referee in the best case scenario, i.e., when the $\varepsilon$-admissibility of each $(x_H, y_H)$ is not revoked.
If there are multiple referees, the computational cost for the benchmarking can exceed that of the optimization process.

\subsection{Referees settings}\label{subsec:Ref_settings}

The refereeing procedure is sensitive to the configuration of the referees.
The choice of the starting point for the lower-level re-optimization critically influences the outcome of the refereeing procedure.
Likewise, the input parameters supplied to the referee directly condition the lower-level optimization process.

\paragraph{Initialization of the starting point:}
A critical design choice in the refereeing procedure concerns the initialization of the referee $r$ when it challenges a point $(x_k, y_k)$ produced by algorithm $a$. 
A first approach is to systematically select the same starting point as the one employed by algorithm $a$ to produce a lower-level response.
In other words, if $a$ is initialized at $y_k^0$ and produces the solution $(x_k, y_k)$, $r$ is initialized at $y_k^0$ as well.
A straightforward way to apply this strategy is to use the previously computed lower-level response, i.e., set $y_k^0 = y_{k-1}$.
To this end, algorithm $a$ needs to track all the initial points it used to generate each lower-level response.
Nevertheless, it enables the referee $r$ to be in similar conditions as the lower-level algorithm employed by $a$ in the original optimization process.

As a second approach, initializing the referee at $y_k$ promotes an intensified local search in the neighborhood of the challenged point.
Notice that the total function evaluation budget is frequently exhausted during the original run in a DFO setting.
Consequently, this approach may revoke a large proportion of challenged points, possibly yielding $\mathcal{H}_{a,p}^{\Refs} = \emptyset$ regardless of the true quality of the benchmarked algorithms.

A third approach can be applied when the problem instance $p \in \Pbank$ supplies a nominal starting point $(x_0, y_0)$, the referee can be initialized at $y_0$ for all challenged points.
This helps mitigate $r$ being entrapped in the same local minima as $a$, allowing for the possibility of computing better alternative lower-level solutions by looking for unexplored feasible regions.

\paragraph{Parameters of referees:}

Settings differences between the referee and the algorithms in $\A$ directly affect the revocation of challenged points, specifically parameters related to stopping conditions.
Allocating a larger evaluation budget to the referee might identify better lower-level solutions at the cost of additional computational expenses.
Conversely, reducing the budget or relaxing termination tolerances lowers the computational overhead of the refereeing procedure but compromises the $\varepsilon$-admissibility verification.

Method-specific parameters can critically affect the outcome of the refereeing procedure, e.g., the mesh shrinkage parameter for mesh-based direct search methods or the initial trust-region radius of trust-region methods.
Adjusting parameters differently between an internal referee $r$ and $a$ might result in unexpected outcomes in the refereeing procedure.
A straightforward choice in the setting of an internal referee is to select all of its parameters to be identical to those selected for $a$. 

The chosen settings for a given referee $r$ need to remain identical for all algorithms $a \in \A$ to ensure a fair refereeing procedure.
Furthermore, other possible refereeing settings may be used in addition to the aforementioned ones to modulate the outcomes of the refereeing procedure.
These settings can also be combined to calibrate the refereeing procedure to user-specific requirements.

\subsection{Computational effort and evaluation budget unit}
\label{subsec:FollowerEndeavour}

{\bl DFO benchmarking techniques measure computational cost by the number of function evaluations.
To extend these techniques to the bilevel setting, an effort metric must count both upper- and lower-level evaluations in order to reflect the computational effort of a BL-DFO algorithm.}

For a fixed upper-level variable, generating its associated lower-level response entails multiple evaluations of the lower-level problem, which constitutes a substantial computational expense in solving~\eqref{eq:MainBilevelProb_Optimistic}.
Besides, upper- and lower-level evaluation costs might differ in a blackbox optimization context.
For these reasons, one can consider scaling both upper- and lower-level evaluations to account for their relative computational costs.

Consider a scaled computational effort that aggregates upper-level and lower-level function evaluations as
\begin{equation}\label{eq:New_N_BLDFO_LLpov}
    N \coloneqq \lambda N_{UL} + N_{LL},
\end{equation}
where \(N_{UL}\) and \(N_{LL}\) represent the numbers of evaluations of the upper-level and the lower-level problems, respectively.
The scaling parameter \(\lambda \ge 0\) represents the relative cost of evaluating the upper-level and the lower-level.
This aggregation is similar to the one proposed in a surrogate-assisted framework~\citep{G-2025-36}.

When the upper-level evaluation cost is insignificant compared that of the lower-level, setting \(\lambda = 0\) enables to account for lower-level evaluations only.
One has $\lambda \ge 1$ if the upper-level is more onerous than the lower-level,  and $\lambda = 1$ implies that evaluating both levels involves the same computational effort, as assumed in~\citep{Ekmekcioglu2024,Islam2017,Islam2018}.
The value of $\lambda$ can be determined using information related to the instances, e.g., by computing the ratio $t_{UL}/t_{LL}$ where $t_{UL}$ and $t_{LL}$ are the elapsed CPU times to evaluate the upper-level and the lower-level, respectively.
Otherwise, $\lambda$ can be fixed by the user to enable adaptability to instance-specific characteristics.

An alternative scaled computational effort can be defined by
\begin{equation}\label{eq:New_N_BLDFO_ULpov}
    N \coloneqq N_{UL} + \lambda^{-1}N_{LL}.
\end{equation}
Both~\eqref{eq:New_N_BLDFO_LLpov} and~\eqref{eq:New_N_BLDFO_ULpov} provide measures of the overall computational cost involved by a bilevel algorithm.
The scaling factor $\lambda$ or $\lambda^{-1}$ is used for indicating whether the upper- or the lower-level problem is the most costly to evaluate.
Formulation~\eqref{eq:New_N_BLDFO_LLpov} is  adapted to situations in which evaluating the upper-level problem is more onerous than evaluating the lower-level.
Conversely, formulation~\eqref{eq:New_N_BLDFO_ULpov} is preferred {\bl when the lower-level functions are more expensive to evaluate than the upper-level ones.} 

The scaled computational effort $N$, with either formulation~\eqref{eq:New_N_BLDFO_LLpov} or~\eqref{eq:New_N_BLDFO_ULpov}, accounts for the efforts required at both levels.
However, the dimensions of both problems remain to be considered by groups of budget units for data profiles to extend them to BLO.
Data profiles constructed by~\citet{DiKuRiZe2024} employ groups of $n_x+1$ evaluations, counting solely upper-level evaluations, which do not properly represent the global effort deployed to solve~\eqref{eq:MainBilevelProb_Optimistic}.
\citet{Cesaroni2026} constructed data profiles using groups of $n_xn_y + 1$ lower-level function evaluations as effort scaling.
This approach accounts for the dimensions of both problems involved in~\eqref{eq:MainBilevelProb_Optimistic}, but does not enable a straightforward generalization of data profiles: for single-level problems, $n_y=0$ resulting in $n_x n_y + 1 = 1$.

Groups of $(n_x+1)(n_y+1)$ budget units are proposed to extend data profile evaluation groups to BLO for two reasons.
The first one is that it naturally generalizes the groups of $n_p+1$ evaluations proposed by~\citet{MoWi2009} for single-level optimization, where $n_p$ is the dimension of problem instance $p$.
Indeed, if there is no lower-level problem (i.e., $n_y=0$), then the groups considered are of $n_x+1$ (upper-level) evaluations, which corresponds exactly to groups of $n_p+1$ evaluations.

The second reason is that $(n_x+1)(n_y+1)$ represents the minimal number of iterations required by a single unsuccessful step in a bilevel direct search poll phase.
In other words, the minimal effort deployed by a BL-DFO algorithm employing a direct search scheme at both levels to terminate the optimization of problem~\eqref{eq:MainBilevelProb_Optimistic} is of $(n_x+1)(n_y+1)$ evaluations.
Such groups, therefore, represent the same minimal budget required by a single-level direct-search algorithm to terminate a poll phase with a minimal positive basis, i.e., $n_p+1$ evaluations.
Otherwise, existing benchmarking techniques remain unchanged, with the filtered data $\mathcal{H}_{a,p}^{\Refs}$ and the scaled effort $N$ used in place of their single-level counterparts.
The refereeing procedure challenges the admissibility of points to discard non-optimal solutions, while the scaled computational effort $N$ more accurately reflects the true computational cost to solve~\eqref{eq:MainBilevelProb_Optimistic}.


\section{Numerical illustrations}
\label{subsec:ComputationalTests}

In this section, the proposed refereeing strategies, as well as the scaled computational effort, are illustrated on a set of bilevel problems in order to assess their impact on existing benchmarking techniques commonly used in DFO.
The set of instances is taken from the BOLIB Matlab library~\citep{Zhou2021}, which consists of analytical bilevel problems. Except for the internal refereeing strategy, in which the identified DFO solvers are explicitly cited, the solvers are anonymized as \texttt{Algo1}, \texttt{Algo2}, and \texttt{Algo3} throughout the tests.

\subsection{Implementation details}\label{subsec:ImplementationDetails}

The refereeing strategies and tested solvers are all implemented in the Julia $1.11$ programming language.
For this work, a Julia variant of the BOLIB~\citep{Zhou2021} library, entitled \texttt{BOLIB.jl}, is developed to make numerous bilevel problems available in open source.
\texttt{BOLIB.jl} contains a majority of BOLIB instances under the \texttt{BilevelProblem} mutable structure.
Tests are conducted on a standard desktop computer with 64 GB of RAM.

The initial point $(x_0, y_0)$ corresponds to the one given in the BOLIB instance definition.
The budgets of function evaluations are $N_{UL}^{tot} = 300$ and $N_{LL}^{tot} = 100n_y$ for the upper-level and the lower-level, respectively.
The tolerance used in the refereeing procedures to revoke the $\varepsilon$-admissibility of points is set as $\varepsilon = 10^{-9}$.
A two-phase extreme barrier scheme~\citep[Chapter 10]{AuHa2026} 
is used in the situations in which the starting point is infeasible.


The {End-Point}, {Complete} and {Reverse} refereeing procedures are compared under specific settings.
The initialization of each challenged point uses the lower-level variable $y_0$ provided by the BOLIB instance.
Each refereeing procedure uses the same evaluation budget and stopping tolerance that were used by the algorithms of $\mathcal A$ for the lower-level problem.
All other algorithmic parameters are set to the default values, 
except when specific setting modifications are explicitly mentioned.


The value of the scaling parameter $\lambda$ is straightforwardly set as
 $   \lambda \coloneqq \frac{t_{UL}}{t_{LL}}$,
where $t_{UL}$ and $t_{LL}$ are the CPU times to call the upper-level and the lower-level functions, respectively.
A different value of $\lambda$ is used in~\Cref{subsec:lambda_impact} to specifically illustrate its impact.
Since the computational expenses to evaluate the upper- and the lower-level problems are equivalent for a majority of BOLIB instances, formulation~\eqref{eq:New_N_BLDFO_LLpov} of the scaled computational effort $N$ is  chosen for all instances.


The accuracy value~\cite{G-2025-36} used to construct the profiles is defined by
 $   F_{acc}^N \coloneqq \frac{F^N - F^0}{F^* - F^0}$,
where $F^N$ is the value of $F$ after a deployed effort of $N$ effort units, $F^0$ is the largest initial $\varepsilon$-admissible value of $F$ and $F^*$ is the best $\varepsilon$-admissible value of $F$ found among all solvers.
The precision values $\tau \in \{10^{-1}, 10^{-2}, 10^{-4}\}$ are used to determine whether an instance is $\tau$-solved (i.e., $F_{acc}^N \ge 1-\tau$) or not. The percentage of $\tau$-solved instances is displayed on the $y$-axis of the data profiles.


\subsection{Impact of scaling the computational effort}\label{subsec:lambda_impact}

The first test set illustrates the impact of employing $N = \lambda N_{UL} + N_{LL}$ against $N = N_{UL}$ and $N= N_{LL}$. 
The context is a situation where evaluating the upper-level is significantly more expensive than evaluating the lower-level:
 the parameter $\lambda$ is set to $100$ for problem instances.
Two anonymous algorithms with different lower-level subsolvers, named \texttt{Algo1} and \texttt{Algo2}, are compared and the results are illustrated with data profiles in~\Cref{fig:data_profile_Compare_Lambda} on $67$ BOLIB instances.
The central plot uses the scaled computational effort $N = \lambda N_{UL} + N_{LL}$ to represent the actual cost spent to solve~\eqref{eq:MainBilevelProb_Optimistic}.
The data profiles on the left considers only the upper-level function evaluations (i.e., $N = N_{UL}$) while the plot on the right considers only the lower-level function evaluations (i.e., $N = N_{LL}$).

\globallegendCompareLambda
\begingroup
\begin{figure}[ht]
    \centering
    \includegraphics[width=1.0\linewidth, trim={5mm 7mm 7mm 5mm},clip]{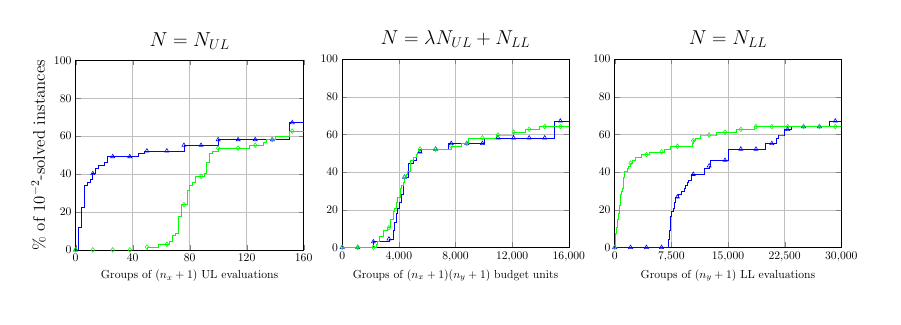}
    \begin{center}
      \pgfplotslegendfromname{globallegendCompareLambda}
    \end{center}
    \caption{{\bl Data profiles using different budget units.}}
    \label{fig:data_profile_Compare_Lambda}
\end{figure}
\endgroup

Notice that the groups of budget units on the x-axis in~\Cref{fig:data_profile_Compare_Lambda} are not identical, resulting in different axis labels.
For $N = N_{UL}$, the effort measured relies only on upper-level evaluations.
In this case, to adapt the x-axis for the left data profile, evaluations are measured in groups of $(n_x+1)$ UL evaluations.
Conversely, with $N = N_{LL}$, the effort measured relies only on the LL evaluations; therefore, groups of $(n_y+1)$ LL evaluations are considered.
The scaled computational effort $N$ is defined by~\eqref{eq:New_N_BLDFO_LLpov} and groups of $(n_x+1)(n_y+1)$ are considered, as detailed in~\Cref{subsec:FollowerEndeavour}.
The budget unit for the data profiles in the central column refers to lower-level equivalent evaluations, as $N$ defined by~\eqref{eq:New_N_BLDFO_LLpov} quantifies the deployed effort from a lower-level perspective.

{\bl When effort is measured by upper-level evaluations only (left plot), \texttt{Algo1} appears to dominate \texttt{Algo2}, whereas the opposite holds when only lower-level evaluations are counted (right plot).
The central plot, which uses the combined effort metric $N = \lambda N_{UL} + N_{LL}$, strikes a balance between these two extremes: under this metric, the data profiles of \texttt{Algo1} and \texttt{Algo2} are equivalent.
Consequently, neglecting either the upper- or the lower-level effort may unduly favor some methods and lead to inaccurate comparisons between benchmarked algorithms.}



\subsection{External referee}\label{subsec:ExternalReferee}


The following example demonstrates that omitting the refereeing procedure can lead to erroneous conclusions.
Three solvers named \texttt{Algo1}, \texttt{Algo2}, and \texttt{Algo3}, differing only in the employed subsolver to address the lower-level problem, are compared.
For each variant, an External Referee is used to solve the lower-level problem  and challenges the solutions generated by employing the End-Point, Reverse, and Complete refereeing strategies.

Results are reported by data profiles in~\Cref{fig:data_profile_AnonymousReferee_LL}. 
The plots in the left column compare the algorithms without any refereeing procedure for two values of $\tau$.
In other words, it is assumed that all solutions produced by the algorithms are $\varepsilon$-admissible.
The three other columns involve a referee that challenges the $\varepsilon$-admissibility of some solutions.
The effect is drastic: the $\varepsilon$-admissibility of many solutions produced by \texttt{Algo1} and \texttt{Algo2} is revoked by the referee.
As a consequence, the data profiles of these two algorithms are lowered, which benefits the performance of \texttt{Algo3}.
With any of the three refereeing processes, \texttt{Algo3} is clearly shown to be dominant. 
\Cref{fig:data_profile_AnonymousReferee_LL} also suggests that {\bl the three refereeing strategies yield similar outcomes, as the second, third and fourth columns of the figure are comparable.
In all cases, removing non-admissible solutions prevents misinterpretations of algorithmic performance.}

\globallegendExternReferee

\begingroup
\begin{figure}[ht]
    \centering
    \includegraphics[width=1.0\linewidth, trim={5mm 7mm 7mm 5mm},clip]{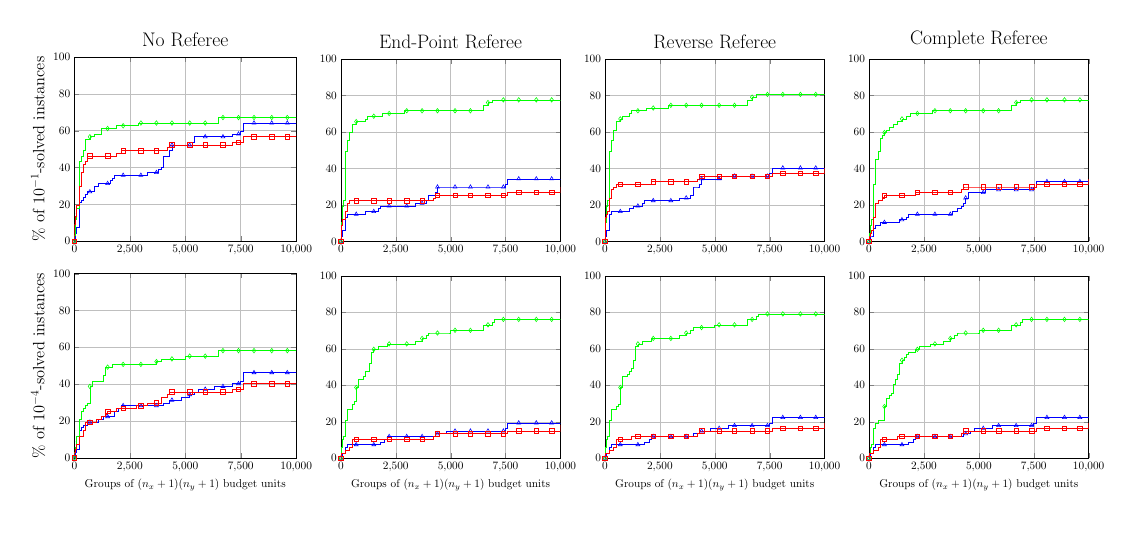}
    \begin{center}
      \pgfplotslegendfromname{globallegendExternReferee}
    \end{center}
    \caption{{\bl Data profiles obtained using an external referee.}}
    \label{fig:data_profile_AnonymousReferee_LL}
\end{figure}
\endgroup

In summary, applying {\bl a refereeing procedure}, even the inexpensive referee that challenges a single point, is necessary to benchmark {\bl BL-DFO} algorithms.
{\bl The refereeing strategy choice depends on the purpose of the refereeing procedure itself.}

\subsection{Internal referee}\label{subsec:InternalReferee}

{\bl This section illustrates the effect of an internal refereeing procedure on three algorithms that share the same upper-level method but use lower-level subsolvers of increasing quality.
\texttt{Algo1} employs a low-quality lower-level subsolver, \texttt{Algo2} uses a more accurate one, and \texttt{Algo3} uses a state-of-the-art one.
\Cref{fig:data_profile_KnownReferee_LL} shows the resulting data profiles.
The first column compares the algorithms without refereeing, i.e., assuming that all the solutions they produce are $\varepsilon$-admissible, while each of the remaining columns applies one refereeing strategy to verify the $\varepsilon$-admissibility of these solutions.
Without refereeing, the three algorithms perform similarly.
Applying any refereeing strategy, however, drastically lowers the data profiles of \texttt{Algo1} and \texttt{Algo2}, whereas \texttt{Algo3} is only marginally affected.
Refereeing thus reveals the hierarchy between the lower-level subsolvers, which is hidden when the reported solutions are taken at face value.

\globallegendInternReferee
\begingroup
\begin{figure}[ht]
    \centering
    \includegraphics[width=1.0\linewidth, trim={5mm 7mm 7mm 5mm},clip]{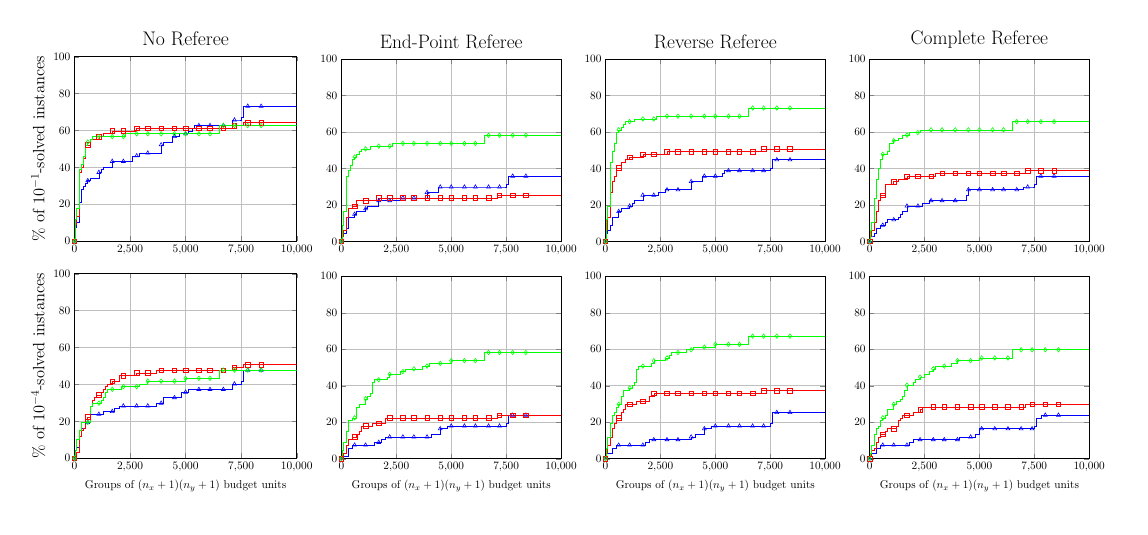}
    \begin{center}
      \pgfplotslegendfromname{globallegendExternReferee}
    \end{center}
    \caption{{\bl Data profiles obtained using an internal referee.}}
    \label{fig:data_profile_KnownReferee_LL}
\end{figure}
\endgroup

Inspection of the logs reveals that the referee, using the lower-level subsolver of \texttt{Algo3}, revoked the $\varepsilon$-admissibility of many points produced by \texttt{Algo1} and \texttt{Algo2}, whose lower-level responses were often not optimal.
Since \texttt{Algo3} employs the most accurate subsolver, its data profiles are much less affected than those of \texttt{Algo1} and \texttt{Algo2}.
Consequently, all three refereeing strategies indicate that \texttt{Algo3} dominates \texttt{Algo1} and \texttt{Algo2}.

The three refereeing strategies nevertheless differ in some respects.
The End-Point referee overestimates the performance of \texttt{Algo1}, whereas the Reverse and Complete referees agree on the same ranking: \texttt{Algo3} dominates \texttt{Algo2}, which in turn outperforms \texttt{Algo1}.
Compared with the Complete referee, the Reverse referee slightly overestimates the performance of all three algorithms, without affecting the general conclusions drawn from the data profiles.}

{\bl This example shows that applying any refereeing strategy changes the conclusions drawn from data profiles.
Moreover, the more expensive refereeing strategies clearly recover the hierarchy of subsolvers stated at the beginning of this section: \texttt{Algo3} outperforms \texttt{Algo2}, which in turn outperforms \texttt{Algo1}.
Although the End-Point referee correctly identifies \texttt{Algo3} as the best algorithm, a more expensive refereeing strategy is recommended for a finer comparison between BL-DFO algorithms.}

{\bl 
\subsection{Computational effort for data collection and for refereeing processes}\label{subsec:PerspectivesTimes}

The CPU times required to solve the BLO problems and to apply all refereeing procedures are displayed in~\Cref{fig:HistogramCPUTimes}. Two features keep the cost of refereeing low. 
First, the referees only perform optimization at feasible points claimed to improve the upper-level objective value. 
Second, because of the intersection in Equation~\eqref{eq:Hap_filtered}, the refereeing process terminates as soon as a referee identifies a non-admissible point. 
As a result, every refereeing procedure requires less time than the original data collection, even when all algorithms take part in the refereeing. 
For both the external and internal referee procedures, the End-Point Referee is by far the cheapest strategy. The Reverse Referee is more expensive, and the Complete Referee is the most expensive of the three.

\begin{figure}[ht]
    \centering
    \begin{subfigure}[b]{0.48\linewidth}
        \centering
        \includegraphics[width=\linewidth]{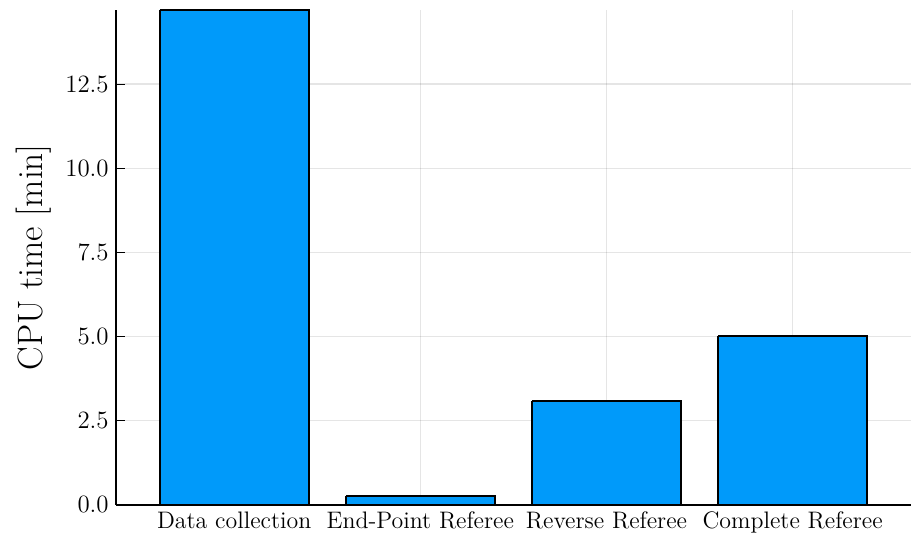}
        \caption{External referee (\cref{subsec:ExternalReferee})}
        \label{fig:HistogramCPUTimesExternal}
    \end{subfigure}
    \hfill
    \begin{subfigure}[b]{0.48\linewidth}
        \centering
        \includegraphics[width=\linewidth]{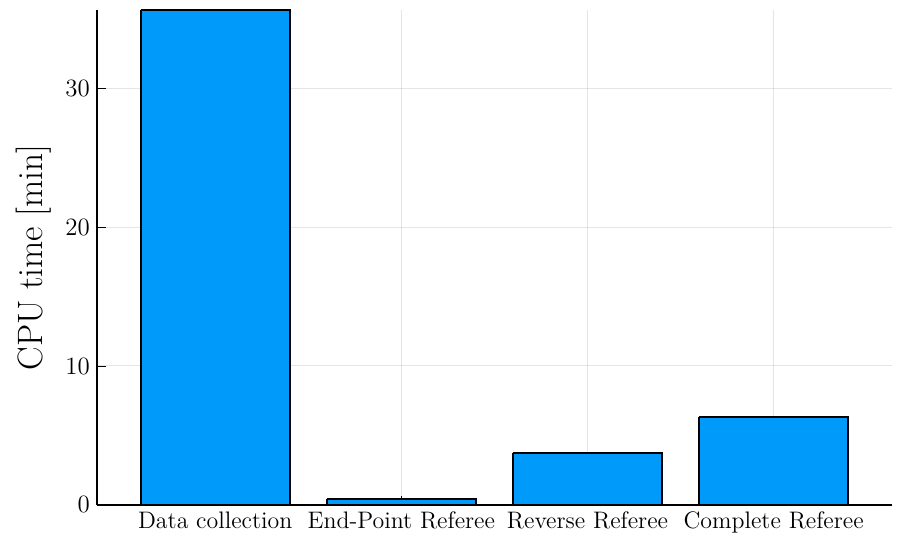}
        \caption{Internal referee (\cref{subsec:InternalReferee})}
        \label{fig:HistogramCPUTimesInternal}
    \end{subfigure}
    \caption{Computational times for data collection and each referee procedure.}
    \label{fig:HistogramCPUTimes}
\end{figure}
}

\section{Conclusion}\label{sec:Conclusion}

This work introduces an extension of existing single-level DFO benchmarking techniques to the bilevel setting. A refereeing procedure is proposed to challenge the $\varepsilon$-admissibility of solutions generated by BL-DFO algorithms, implemented through three strategies of increasing computational cost:
{\bl the End-Point referee challenges only the final point of each run and is thus the least expensive but also the least accurate; the Complete Referee challenges every point in the log and is the most reliable, at the highest computational cost; and the Reverse referee offers a trade-off between the two, revoking points from the end of the log up to a specific threshold.}

A scaled computational effort aggregating both upper- and lower-level evaluations is also introduced to better represent the overall effort required to solve nonlinear bilevel problems. Both contributions are illustrated on analytical test problems. The refereeing procedure is shown to discard non-admissible points, yielding more faithful conclusions regarding algorithmic performance.
Similarly, the scaled computational effort provides a more accurate measure of the total effort deployed, as neglecting part of this effort may lead to overestimated performances in data profiles.



\small
  \bibliographystyle{abbrvnat.bst}
  \bibliography{bibliography}
\normalsize

\end{document}